\newcommand{\doublespace} {\addtolength{\baselineskip}{.50\baselineskip}}

\documentclass[reqno,10pt]{amsart}
\usepackage {amsmath, amssymb, enumerate, amsthm}
\usepackage{graphicx}
\usepackage{color}
\usepackage{cite}

\newtheorem{theorem}{Theorem}[section]

\newtheorem{definition}{Definition}[section]

\newtheorem{thm}[theorem]{\bf Theorem}

\newtheorem{defn}[definition]{\bf Definition}

\begin{document}

\allowdisplaybreaks  

\doublespace

\title[wave equation with logarithmic source] {Stabilization for the wave equation with fully subciritical logarithmic nonlinearity}

\author[Tae Gab Ha]{Tae Gab Ha}

\address{Department of Mathematics and Institute of Pure and Applied Mathematics, Jeonbuk National University, Jeonju 54896, Republic of Korea}

\email{tgha@jbnu.ac.kr}




\date{}
\maketitle

\begin{abstract}
In this paper, we consider a wave equation with strong damping and logarithmic nonlinearity. This paper aims to study the local and global existence, uniqueness and the uniform energy decay rate of a weak solution under some sufficient conditions on the initial data. Unlike previous literature restricted to the lower subcritical range $2 < \gamma < \frac{2(n-1)}{n-2}$, we successfully extend the validity of the well-posedness and stabilization results to the upper subcritical range $\frac{2(n-1)}{n-2} \leq \gamma < \frac{2n}{n-2}$.

\end{abstract}

\section{ Introduction }
\setcounter{equation}{0}

In this paper, we investigate the local and global existence, uniqueness and energy decay rate of the solution to the following wave equation with logarithmic nonlinearity:
\begin{equation}\label{1}
\begin{cases}
\vspace{3mm} u_{tt} -  \Delta u - \Delta u_t = |u|^{\gamma -2} u \ln|u| &\hspace{5mm} \text{in} \hspace{5mm} \Omega~~\times~~ (0, +\infty),
\\
\vspace{3mm} u = 0 &\hspace{5mm} \text{on} \hspace{5mm} \partial\Omega ~~\times~~ (0, +\infty),
\\
\vspace{3mm} u(x,0) = u_0(x), \hspace{5mm} u_t(x,0) = u_1(x) &\hspace{5mm} \text{in} \hspace{5mm} \Omega,
\end{cases}
\end{equation}
where $\Omega \subset \mathbb{R}^n$ $(n \geq 3)$ is a bounded domain with smooth boundary $\partial \Omega$, and $2 < \gamma < \frac{2n}{n-2}$.

Nonlinear evolution equations, such as the wave equations, provide fundamental models for describing a wide array of phenomena in physics, chemistry, and engineering. The inclusion of nonlinear source terms in these equations leads to a rich and complex dynamic, often characterized by a competition between the natural dissipative or dispersive effects of the system and the concentrating effect of the nonlinearity. This competition can result in solutions that exist globally in time or, conversely, develop singularities and blow up in finite time.

A particularly interesting class of nonlinearities that has garnered significant attention is the logarithmic source term. This type of nonlinearity appears in various physical contexts, including quantum mechanics and nuclear physics (cf. \cite{biru1, biru2, heft} ). In recent years, significant attention has been devoted to wave equations with logarithmic nonlinearities of the form $|u|^{\gamma -2} u \ln|u|$ (see \cite{di, guo, ha, wu} and reference therein).  For instance, Guo and Zhang \cite{guo} studied general decay results for a viscoelastic wave equation with the logarithmic nonlinear source and dynamic Wentzell boundary condition. Di et al \cite{di} analyzed \eqref{1}, establishing existence and blow-up results under certain conditions. However, all references mentioned above have restricted their focus to exponents in the lower subcritical range, $2 < \gamma < \frac{2(n-1)}{n-2}$. For exponents beyond this classical threshold, standard energy estimates and classical uniqueness proofs break down due to a critical loss of compactness and the inability to bound the nonlinear difference in standard $L^2$ spaces. Consequently, the well-posedness of the wave equation with a logarithmic source in the upper subcritical range $\frac{2(n-1)}{n-2} \leq \gamma < \frac{2n}{n-2}$ has remained an open problem.

This paper aims to fill this gap by studying the wave equation with a logarithmic source term over the fully subcritical range up to the Sobolev critical exponent.
This paper show that in the upper subcritical range, $\frac{2(n-1)}{n-2} \leq \gamma < \frac{2n}{n-2}$ the same results can be obtained as in the lower subcritical range. The main objective of this paper is to establish the local and global existence, uniqueness, and uniform energy decay rate for weak solutions to this problem.

The paper is organized as follows. Section 2 introduces the notation, functional setting, and main theorem. Section 3 is dedicated to the proof of the main theorem, broken down into subsections covering local existence, global existence, uniqueness, and the energy decay estimate.

\section{ Preliminaries }
\setcounter{equation}{0}

We begin this section by introducing some notations and definitions that will be used throughout the paper. Let $\Omega \subset \mathbb{R}^n$ $(n \geq 3)$ be a bounded domain with smooth boundary $\partial\Omega$. Throughout the paper we use standard functional spaces and denote that $\|\cdot\|_p$ is the $L^p(\Omega)$-norm and $(u,v) = \int_\Omega u(x) v(x) dx$.

\textbf{Assumption (A)}. The exponent $\gamma$ satisfies
\begin{equation}\label{201}
\frac{2(n-1)}{n-2} \leq \gamma < \frac{2n}{n-2}.
\end{equation}

The energy associated to problem \eqref{1} is given by
\begin{equation}\label{202}
E(t) := E(u,u_t) = \frac{1}{2}\|u_t\|^2_2 + \frac{1}{2}\|\nabla u\|^2_2 - \frac{1}{\gamma}\int_\Omega |u|^\gamma \ln|u| dx + \frac{1}{\gamma^2}\|u\|^\gamma_\gamma.
\end{equation}
To deal with logarithmic source term, we define some functional as follows:
\begin{equation}\label{203}
J(u) =  \frac{1}{2}\|\nabla u\|^2_2 - \frac{1}{\gamma}\int_\Omega |u|^\gamma \ln|u| dx + \frac{1}{\gamma^2}\|u\|^\gamma_\gamma,
\end{equation}
\begin{equation}\label{204}
I(u) = \|\nabla u\|^2_2 - \int_\Omega |u|^\gamma \ln|u| dx,
\end{equation}
\begin{equation}\label{205}
d = \inf_{u \in H^1_0(\Omega)\setminus\{0\}} \sup_{\lambda >0} J(\lambda u) = \inf_{u \in \mathcal{N}} J(u), ~~\text{where}~~ \mathcal{N} = \{u \in H^1_0(\Omega)\setminus \{0\} ~|~ I(u) = 0\}.
\end{equation}
Here $\mathcal{N}$ is called the Nehari manifold, and $d$ is the depth of potential well. It is well known that $\mathcal{N}$ is not empty and $d$ is positive (see \cite{chen2, di}).



\begin{defn} (Weak solution)
The function $u$ is called a weak solution of problem \eqref{1} on $\Omega \times (0,T)$ if $u \in C(0,T; H^1_0(\Omega)) \cap C^1 (0, T ; L^2(\Omega)) \cap C^2(0,T ; H^{-1}(\Omega))$ with $u_t \in L^2(0,T ; H^1_0(\Omega))$ and satisfies problem \eqref{1} in the distribution sense, i.e.
$$
(u_{tt}, \phi) + (\nabla u, \nabla \phi) + (\nabla u_t, \nabla \phi) = (|u|^{\gamma -2} u \ln |u|, \phi),
$$
for any $\phi \in H^1_0(\Omega)$, and $u(x,0) = u_0$ in $H^1_0(\Omega)$, $u_t(x,0) = u_1$ in $L^2(\Omega)$.
\end{defn}

\begin{defn} (Maximal existence time)
Let $u$ be a weak solution of problem \eqref{1} on $\Omega \times (0,T)$, we define the maximal existence time $T_{\max}$ as follows
$$
T_{\max} = \sup \{T > 0 ; ~u ~~\text{exists on}~~ (0,T)\}.
$$
\begin{enumerate}[(i)]
\item If $T_{\max} = +\infty$, we say that the solution $u$ is global.
\item If $T_{\max} < +\infty$, we say that the solution $u$ blows up in finite time and $T_{\max}$ is the blow-up time.
\end{enumerate}
\end{defn}

We now state our main result.

\begin{thm}
Let $u_0 \in H^1_0(\Omega)$, $u_1 \in L^2(\Omega)$, $I(u_0) > 0$ and $E(0) < d$. Assume that (A) is hold, then the problem \eqref{1} admits a unique global weak solution. Furthermore, there exist positive constants $C_1$ and $C_2$ such that $E(t)$ satisfies
$$
E(t) \leq C_1 e^{-C_2 t} ~~\text{for all}~~ t >0.
$$
\end{thm}

\section{ Proof of Theorem 2.1 }
\setcounter{equation}{0}

\subsection{ Local existence }

Our goal in this subsection is to show the local existence of weak solutions to problem \eqref{1} using the Faedo-Galerkin method. Let $\{w_j\}_{j \in \mathbb{N}}$ be an orthogonal basis of $H^1_0(\Omega)$ and define $V_m = span\{w_1, w_2, \cdots, w_m \}$. We search for an approximate solution, for each $m \in \mathbb{N}$,
$$
u^m(t) = \sum^m_{j=1} \delta^m_j(t)w_j
$$
satisfying the approximate equation

\begin{equation}\label{301}
\begin{cases}
\vspace{3mm} (u^m_{tt}, \varphi) + (\nabla u^m, \nabla \varphi) + (\nabla u^m_t, \nabla \varphi) = (|u^m|^{\gamma -2} u^m \ln |u^m|, \varphi), \hspace{3mm} \varphi \in V_m,
\\
\vspace{3mm} u^m(0) = u^m_0 = \sum^m_{j=1} (u_0, w_j)w_j \rightarrow u_0 \hspace{3mm}\text{strongly in}\hspace{3mm} H^1_0(\Omega),
\\
 u^m_t(0) = u^m_1 = \sum^m_{j=1} (u_1, w_j)w_j \rightarrow u_1 \hspace{3mm}\text{strongly in}\hspace{3mm} L^2(\Omega).
\end{cases}
\end{equation}
Theory of ordinary differential equations provides a unique local solution $u^m$ to problem \eqref{301} on the internal $[0, t_m)$, $t_m \in (0, T]$. A solution $u$ to problem \eqref{1} on some internal $[0, t_m)$ will be obtained as the limit of $u^m$ as $m \rightarrow \infty$. Next, we show that $t_m = T$ and the local solution is uniformly bounded independent of $m$ and $t$. For this purpose, let us replace $\varphi$ by $u^m_t$ in \eqref{301} we obtain
\begin{equation}\label{302}
E(u^m, u^m_t) + \int^t_0 \|\nabla u^m_s\|^2_2 ds = E(u^m_0, u^m_1),
\end{equation}
where
$$
E(u^m,u^m_t) = \frac{1}{2}\|u^m_t\|^2_2 + \frac{1}{2}\|\nabla u^m\|^2_2 - \frac{1}{\gamma}\int_\Omega |u^m|^\gamma \ln|u^m| dx + \frac{1}{\gamma^2}\|u^m\|^\gamma_\gamma.
$$
By the continuity of $E(t)$ and $I(t)$, it follows that
\begin{equation}\label{303}
E(u^m_0, u^m_1) < d \hspace{3mm}\text{and}\hspace{3mm} I(u^m_0)>0,
\end{equation}
for sufficiently large $m$. This asserts
\begin{equation}\label{304}
E(u^m(t), u^m_t(t)) < d \hspace{3mm}\text{and}\hspace{3mm} I(u^m(t))>0,
\end{equation}
for sufficiently large $m$ and $t \geq 0$. Indeed, from \eqref{302} and \eqref{303}, it is clear that $E(u^m(t), u^m_t(t)) < d$. On the other hand, if $I(u^m(t))>0$ is false, then there exists $t_0 >0$ such that $I(u^m(t_0))=0$. Then $E(u^m(t_0), u^m_t(t_0)) \geq J(u^m(t_0)) \geq d$, which is contradiction. Therefore from \eqref{302} and \eqref{304} we conclude that
\begin{equation}\label{305}
\frac{1}{2}\|u^m_t\|^2_2 + \frac{\gamma - 2}{2\gamma}\|\nabla u^m\|^2_2 + \frac{1}{\gamma^2}\|u^m\|^\gamma_\gamma + \int^t_0 \|\nabla u^m_s\|^2_2 ds <d.
\end{equation}
Hence \eqref{305} permit us to obtain a subsequence of $\{u^m\}$, which we still denote it by $\{u^m\}$ such that
$$
u^m \rightarrow u  \hspace{3mm}\text{weak star in}\hspace{3mm} L^\infty(0,T ; H^1_0(\Omega)),
$$
$$
u^m_t\rightarrow u_t  \hspace{3mm}\text{weak star in}\hspace{3mm} L^\infty(0,T ; L^2(\Omega)).
$$
$$
u^m\rightarrow u  \hspace{3mm}\text{weak star in}\hspace{3mm} L^\infty(0,T ; L^\gamma(\Omega)).
$$
$$
u^m_t \rightarrow u_t  \hspace{3mm}\text{weakly in}\hspace{3mm} L^2(0,T ; H^1_0(\Omega)),
$$

By Aubin-Lions compactness result (see \cite{simon}, Corollary 4), we have
$$
u^m \rightarrow u  \hspace{3mm}\text{strongly in}\hspace{3mm} C(0,T ; L^2(\Omega)),
$$
which implies that
$$
u^m \rightarrow u  \hspace{3mm}\text{a.e.}\hspace{3mm} \Omega \times (0,T)
$$
and
\begin{equation}\label{306}
|u^m|^{\gamma -2} u^m \ln |u^m| \rightarrow |u|^{\gamma -2} u \ln |u|  \hspace{3mm}\text{a.e.}\hspace{3mm} \Omega \times (0,T).
\end{equation}
Next, we show that
\begin{equation}\label{307}
|u^m|^{\gamma -2} u^m \ln |u^m|  \hspace{3mm}\text{is bounded in}\hspace{3mm} L^\infty(0,T ; L^\frac{\gamma}{\gamma - 1}(\Omega)).
\end{equation}
For this, we let
$$
\Omega_1 = \{x \in \Omega ~:~ |u^m| < 1\}  \hspace{3mm}\text{and}\hspace{3mm} \Omega_2 = \{x \in \Omega ~:~ |u^m| \geq 1\}.
$$
Since $\frac{2(n-1)}{n-2} \leq \gamma < \frac{2n}{n-2}$, there exists a positive constant $\rho$ such that $\rho < \frac{2n}{n-2} - \gamma$. We set a positive constant $\mu = \frac{\rho(\gamma-1)}{\gamma}$. Then we have $\frac{2(n-1)}{n-2} \leq \gamma + \frac{\mu\gamma}{\gamma - 1} < \frac{2n}{n-2}$. From the fact $|s^{\gamma -1} \ln s| \leq (e(\gamma - 1))^{-1}$ for any $0 < s < 1$, while $s^{-\mu} \ln s <  (e \mu)^{-1}$ for any $s \geq 1$, $\mu >0$, we obtain
\begin{equation}\label{308}
\begin{aligned}
&\int_\Omega \Bigl| |u^m|^{\gamma -2} u^m \ln |u^m| \Bigr|^\frac{\gamma}{\gamma-1} dx
\\
&= \int_{\Omega_1} \Bigl| |u^m|^{\gamma -2} u^m \ln |u^m| \Bigr|^\frac{\gamma}{\gamma-1} dx + \int_{\Omega_2} \Bigl| |u^m|^{\gamma -2} u^m \ln |u^m| \Bigr|^\frac{\gamma}{\gamma-1} dx
\\
& \leq \bigl(e(\gamma -1)\bigr)^{-\frac{\gamma}{\gamma-1}}|\Omega| + (e\mu)^{-\frac{\gamma}{\gamma-1}}C^{\gamma + \frac{\mu\gamma}{\gamma - 1}}_* \|\nabla u^m\|^{\gamma + \frac{\mu\gamma}{\gamma - 1}}_2
\\
& \leq \bigl(e(\gamma -1)\bigr)^{-\frac{\gamma}{\gamma-1}}|\Omega| +  (e\mu)^{-\frac{\gamma}{\gamma-1}}C^{\gamma + \frac{\mu\gamma}{\gamma - 1}}_* \Bigl(\frac{2\gamma d}{\gamma -2} \Bigr)^{\frac{1}{2}\bigl(\gamma + \frac{\mu\gamma}{\gamma - 1}\bigr)},
\end{aligned}
\end{equation}
where $C_*$ is the best constant of the Sobolev embedding $H^1_0(\Omega) \hookrightarrow L^{\gamma + \frac{\mu\gamma}{\gamma - 1}}(\Omega)$. Hence by \eqref{306}, \eqref{307} and Lions Lemma (see \cite{lions}, Lemma 1.3), one gets
$$
|u^m|^{\gamma -2} u^m \ln |u^m| \rightarrow |u|^{\gamma -2} u \ln |u| \hspace{3mm}\text{weak star in}\hspace{3mm} L^\infty(0,T ; L^\frac{\gamma}{\gamma-1}(\Omega)).
$$
Therefore, passing to the limit in \eqref{301}, it follows that
$$
(u_{tt}, \varphi) + (\nabla u, \nabla \varphi) + (\nabla u_t, \nabla \varphi) = (|u|^{\gamma -2} u \ln |u|, \varphi),
$$
for all $\varphi \in H^1_0(\Omega)$, and $u(x,0) = u_0$ in $H^1_0(\Omega)$, $u_t(x,0) = u_1$ in $L^2(\Omega)$.

\subsection{ Global existence }

In order to prove the global existence of solutions to problem \eqref{1}, it suffices to show that $\|u_t\|^2_2 + \|\nabla u\|^2_2 + \|u\|^\gamma_\gamma$ is bounded independent of $t$. Multiplying \eqref{1} by $u_t$ and integrating it over $\Omega \times [0,t)$, we have
\begin{equation}\label{309}
E(t) + \int^t_0 \|\nabla u_s\|^2_2 ds = E(0)
\end{equation}
consequently, $E'(t) = -\|\nabla u_t\|^2_2 \leq 0$, which implies that $E(t)$ is a nonincreasing function with respect to $t$. By similar arguments to \eqref{304}, we obtain
\begin{equation}\label{310}
E(t) < d \hspace{3mm}\text{and}\hspace{3mm} I(t) > 0,
\end{equation}
for all $t$. Therefore, by definition of $E(t)$ and $I(t)$, and \eqref{310}, we conclude that
\begin{equation}\label{311}
d > E(0) \geq E(t) > \frac{1}{2} \|u_t\|^2_2 + \frac{\gamma - 2}{2\gamma} \|\nabla u\|^2_2 + \frac{1}{\gamma^2}\|u\|^\gamma_\gamma \geq C_3 \bigl(\|u_t\|^2_2 + \|\nabla u\|^2_2 + \|u\|^\gamma_\gamma \bigr),
\end{equation}
where $C_3 = \min\bigl\{\frac{1}{2}, \frac{\gamma-2}{2\gamma}, \frac{1}{\gamma^2}\bigr\}$. This is the completion of the proof of the global existence of solutions to problem \eqref{1}.

\subsection{ Uniqueness }

Let $v$ and $w$ be two solutions of problem \eqref{301} which have the same initial data. Then $z = v - w$ verifies
$$
(z_{tt}, \varphi) + (\nabla z, \nabla \varphi) + (\nabla z_t, \nabla \varphi) = (|v|^{\gamma -2} v \ln |v| - |w|^{\gamma -2} w \ln |w| , \varphi),
$$
for all $\varphi \in H^1_0(\Omega)$, and $z(x,0) = z_t(x,0) = 0$. By replacing $\varphi = z_t$ in above identity, we have
\begin{equation}\label{312}
\begin{aligned}
&\frac{d}{dt} \bigl[\|z_t\|^2_2 + \|\nabla z\|^2_2 \bigl] + 2\|\nabla z_t\|^2_2
\\
&= 2\int_\Omega (|v|^{\gamma-2} v \ln |v| - |w|^{\gamma -2} w \ln |w|) z_t dx
\\
& = 2\int_\Omega \bigl((\gamma - 1) |\zeta|^{\gamma-2} \ln |\zeta| + |\zeta|^{\gamma -2} \bigr)z~ z_t dx
\\
& = 2\int_\Omega  |\zeta|^{\gamma-2} z~z_t dx + 2(\gamma -1) \int_\Omega  |\zeta|^{\gamma-2} \ln  |\zeta|~ z~z_t dx
\\& := I_1 + I_2,
\end{aligned}
\end{equation}
where $\zeta = \theta v + (1-\theta)w$, $0< \theta < 1$. From the H\"{o}lder inequality with $\frac{2}{n} + \frac{n-2}{2n} + \frac{n-2}{2n} = 1$ and the Young inequality, we deduce that
\begin{equation}\label{313}
\begin{aligned}
I_1 &= 2\int_\Omega  |\zeta|^{\gamma-2} z~z_t dx
\\
&\leq 2\|\zeta\|^{\gamma-2}_{\frac{n}{2}(\gamma-2)} \|z\|_\frac{2n}{n-2} \|z_t\|_\frac{2n}{n-2}
\\
&\leq C_4 \|\nabla \zeta\|^{\gamma-2}_2 \|\nabla z\|_2 \|\nabla z_t\|_2
\\
& \leq C_{5}(\eta) \|\nabla z\|^2_2 + \eta \|\nabla z_t\|^2_2,
\end{aligned}
\end{equation}
where where $C_{4}$ is a positive constant and $C_{5}(\eta)$ depends on $\eta$. The Sobolev embedding $H^1_0(\Omega) \hookrightarrow L^{\frac{n}{2}(\gamma-2)}(\Omega)$ could be used here because the condition $\frac{n}{2}(\gamma-2) < \frac{2n}{n-2}$ is satisfied. Using again the H\"{o}lder inequality with $\frac{2}{n} + \frac{n-2}{2n} + \frac{n-2}{2n} = 1$, we have
\begin{equation}\label{314}
\begin{aligned}
I_2 &= 2(\gamma -1) \int_\Omega  |\zeta|^{\gamma-2} \ln  |\zeta|~ z~z_t dx
\\
& \leq C_6 \|~ |\zeta|^{\gamma-2} \ln  |\zeta|~\|_\frac{n}{2} \|\nabla z\|_2 \|\nabla z_t\|_2,
\end{aligned}
\end{equation}
where $C_6$ is a positive constant. Since $\frac{2(n-1)}{n-2} \leq \gamma < \frac{2n}{n-2}$, we can choose $\mu_1 >0$ such that $\frac{n}{2} (\gamma - 2 + \mu_1) < \frac{2n}{n-2}$. By similar arguments to \eqref{308}, it follows that
\begin{equation}\label{315}
\int_\Omega |~ |\zeta|^{\gamma-2} \ln  |\zeta|~|^\frac{n}{2} dx \leq \bigl(e(\gamma -2)\bigr)^{-\frac{n}{2}}|\Omega| + (e\mu_1)^{-\frac{n}{2}} C^{\frac{n}{2} (\gamma - 2 + \mu_1)}_{**} (2C^{-1}_3 d)^{\frac{n}{4} (\gamma - 2 + \mu_1)},
\end{equation}
where $C_{**}$ is the best constant of the Sobolev embedding $H^1_0(\Omega) \hookrightarrow L^{\frac{n}{2} (\gamma - 2 + \mu_1)}(\Omega)$. Inserting \eqref{315} into \eqref{314} and applying Young's inequality with the same $\eta>0$, we have
\begin{equation}\label{316}
I_2 \leq C_7(\eta) \|\nabla z\|^2_2 +  \eta\|\nabla z_t\|^2_2,
\end{equation}
where $C_7(\eta)$ is a positive constant. Combining \eqref{312}, \eqref{313} and \eqref{316}, we get
$$
\frac{d}{dt} \bigl[\|z_t\|^2_2 + \|\nabla z\|^2_2 \bigl] + 2(1 - \eta)\|\nabla z_t\|^2_2 \leq (C_5(\eta) + C_7(\eta)) \|\nabla z\|^2_2.
$$
By choosing $\eta>0$ sufficiently small such that $0<\eta<1$, we obtain
\begin{equation*}
\frac{d}{dt}\left[\|z_{t}\|_{2}^{2}+\|\nabla z\|_{2}^{2}\right]\le C_{8}\|\nabla z\|_{2}^{2} \le C_{8}\left(\|z_{t}\|_{2}^{2}+\|\nabla z\|_{2}^{2}\right),
\end{equation*}
where $C_8$ is a positive constant. Integrating it over $[0,t)$ and employing the Gronwall Lemma, we conclude that $z = 0$.

\subsection{ Energy decay }

In this subsection, we prove the uniform exponential energy decay rate for problem \eqref{1}. In the following, $C$ denotes a generic positive constant.

First, integrating the energy dissipation law $E'(t) = -\|\nabla u_{t}\|_{2}^{2} \le 0$ over $[S, T]$, we obtain
\begin{equation}\label{eq3.17}
\int_{S}^{T} \|\nabla u_{t}\|_{2}^{2} dt = E(S) - E(T) \le E(S),
\end{equation}
consequently, by applying the Poincar\'e inequality,
\begin{equation}\label{eq3.18}
\int_{S}^{T} \|u_{t}\|_{2}^{2} dt \le C_{P}^{2} \int_{S}^{T} \|\nabla u_{t}\|_{2}^{2} dt \le C_{P}^{2} E(S),
\end{equation}
where $C_P$ is the Poincar\'e constant.

Next, multiplying equation \eqref{1} by $u$ and integrating it over $\Omega\times[S,T]$ and recalling the definition of $I(u) = \|\nabla u\|_{2}^{2} - \int_{\Omega}|u|^{\gamma}\ln|u|dx$, we have
\begin{equation}\label{eq3.19}
\int_{S}^{T} I(u(t)) dt = \int_{S}^{T} \|u_{t}\|_{2}^{2} dt - \left[ \int_{\Omega} u_{t}u \,dx + \frac{1}{2}\|\nabla u\|_{2}^{2} \right]_{S}^{T}.
\end{equation}

Now, we estimate the second term on the right-hand side of \eqref{eq3.19}. Using the Young and Poincar\'e inequalities, along with the fact from \eqref{311} that $E(t) \ge C_{3}(\|u_{t}\|_{2}^{2} + \|\nabla u\|_{2}^{2})$, we get:
\begin{equation*}
\left| \int_{\Omega} u_{t}u \,dx + \frac{1}{2}\|\nabla u\|_{2}^{2} \right| \le \frac{1}{2}\|u_{t}\|_{2}^{2} + \frac{1+C_{P}^{2}}{2}\|\nabla u\|_{2}^{2} \le C E(t).
\end{equation*}
Since $E(t)$ is non-increasing, we deduce that
\begin{equation}\label{eq3.20}
\left| \left[ \int_{\Omega} u_{t}u \,dx + \frac{1}{2}\|\nabla u\|_{2}^{2} \right]_{S}^{T} \right| \le C E(S) + C E(T) \le 2C E(S).
\end{equation}

Inserting \eqref{eq3.18} and \eqref{eq3.20} into \eqref{eq3.19}, we conclude that
\begin{equation}\label{eq3.21}
\int_{S}^{T} I(u(t)) dt \leq C E(S).
\end{equation}

To complete the proof, we establish an upper bound for the energy $E(t)$ in terms of $\|u_{t}\|_{2}^{2}$ and $I(u(t))$. Using the definition of $I(u(t))$, the total energy can be rewritten as:
\begin{align}\label{eq3.22}
E(t) &= \frac{1}{2}\|u_{t}\|_{2}^{2} + \frac{1}{2}\|\nabla u\|_{2}^{2} - \frac{1}{\gamma}\int_{\Omega}|u|^{\gamma}\ln|u|dx + \frac{1}{\gamma^{2}}\|u\|_{\gamma}^{\gamma} \nonumber \\
&= \frac{1}{2}\|u_{t}\|_{2}^{2} + \frac{\gamma-2}{2\gamma}\|\nabla u\|_{2}^{2} + \frac{1}{\gamma}I(u(t)) + \frac{1}{\gamma^{2}}\|u\|_{\gamma}^{\gamma}.
\end{align}

Since the initial energy satisfies $E(0) < d$ and $I(u_{0}) > 0$, the potential well framework guarantees that the solution remains in the stable set. Thus, there exists a constant $\delta \in (0, 1)$ such that $I(u(t)) \ge \delta \|\nabla u\|_{2}^{2}$ for all $t \ge 0$, which implies:
\begin{equation}\label{eq3.23}
\|\nabla u\|_{2}^{2} \le \frac{1}{\delta} I(u(t)).
\end{equation}

Furthermore, by the Sobolev embedding $H_{0}^{1}(\Omega) \hookrightarrow L^{\gamma}(\Omega)$ and the inequality $\|\nabla u\|_{2}^{2} \le \frac{2\gamma d}{\gamma-2}$ from \eqref{311}, we can bound the $L^{\gamma}$-norm term as follows:
\begin{equation}\label{eq3.24}
\|u\|_{\gamma}^{\gamma} \le C_{*}^{\gamma} \|\nabla u\|_{2}^{\gamma} = C_{*}^{\gamma} \|\nabla u\|_{2}^{\gamma-2} \|\nabla u\|_{2}^{2} \le C_{*}^{\gamma} \left( \frac{2\gamma d}{\gamma-2} \right)^{\frac{\gamma-2}{2}} \|\nabla u\|_{2}^{2} := C_{W} \|\nabla u\|_{2}^{2},
\end{equation}
where $C_{*}$ is the embedding constant and $C_{W} > 0$. Substituting \eqref{eq3.23} into \eqref{eq3.24}, we get $\|u\|_{\gamma}^{\gamma} \le \frac{C_{W}}{\delta} I(u(t))$.

Applying these bounds to \eqref{eq3.22}, we successfully deduce the required upper bound:
\begin{equation*}
E(t) \le \frac{1}{2}\|u_{t}\|_{2}^{2} + \left( \frac{\gamma-2}{2\gamma\delta} + \frac{1}{\gamma} + \frac{C_{W}}{\gamma^{2}\delta} \right) I(u(t)) \le M \left( \|u_{t}\|_{2}^{2} + I(u(t)) \right), \quad \forall t \ge 0,
\end{equation*}
for a sufficiently large constant $M > 0$.

Integrating this inequality over $[S, T]$ and applying \eqref{eq3.18} and \eqref{eq3.21}, we obtain:
\begin{equation}\label{eq3.25}
\int_{S}^{T} E(t) dt \le M \int_{S}^{T} \left( \|u_{t}\|_{2}^{2} + I(u(t)) \right) dt \le M(C_{P}^{2} + C_{4}) E(S) := C_{0} E(S),
\end{equation}
for all $T \ge S \ge 0$, where $C_{0} > 0$ is a constant independent of $S$ and $T$. This implies that $E(t)$ decays to zero exponentially.

\section*{Acknowledgments}

This paper was supported by research funds of Jeonbuk National University in 2024. This research was supported by Basic Science Research Program through the National Research Foundation of Korea(NRF) funded by the Ministry of Education (RS-2022-NR075641).

\end{document}